%% file: artjul01.tex
\documentclass[12pt,a4paper]{article}
\usepackage[latin1]{inputenc}
\usepackage{amssymb}
\usepackage{pictex}
\usepackage{epsfig}

\newtheorem{prop}{Proposition}[section]
\newtheorem{lemma}{Lemma}[section]

\newcommand{\bprf}{{\it Proof.~}}
\newcommand{\eprf}{\hfill $\square$ \bigskip\par}

\newcommand{\pitr}{ \mathbb{P}_3}

\newcommand{\piu}{\mathbb{P}_1}

\newcommand{\ci}{ \mathbb{C}}

\newcommand{\zi}{\mathbb{Z}}
\newcommand{\cix}{\ci[x_0,x_1,x_2,x_3]}
\begin{document}

\include{artjul01_00}


\cleardoublepage

\end{document}

%% file: artjul01_00.tex
\begin{center}
{\Large Polynomial Invariants of the Heisenberg Group\\ with extra Symmetries}\\[0.3cm]
{\rm \today}\\[0.3cm]
 A.Sarti\\[0.3cm]
FBR 17 Mathematik, Universit\"at Mainz\\
Staudingerweg 9, 55099 Mainz\\
Germany\\
e-mail:sarti@mathematik.uni-mainz.de\\
\end{center}

\begin{abstract}
Let $G\subset SO(4)$ denote a finite subgroup containing the Heisenberg group. In these notes we classify all these groups, we find the dimension of the spaces of $G$-invariant polynomials and we give equations for the generators whenever the space has dimension two.  Then we complete the study of the corresponding $G$-invariant pencils of surfaces in $\mathbb{P}_3$ which we started in \cite{sa}. It turns out that we have five more pencils, two of them containing surfaces with nodes.
\end{abstract}

\tableofcontents

\section{Introduction}
Consider the Klein four group $V\subseteq SO(3)$. Let $\tilde{V}$ denotes its inverse image in $SU(2)$ under the universal covering $SU(2)$$\rightarrow$$SO(3)$. The image of the direct product $\tilde{V}\times\tilde{V}$ in $SO(4)$ under the double covering $SU(2)$$\times$$SU(2)$$\rightarrow$$SO(4)$ is the Heisenberg group $H$. In this note we classify all the subgroups $G$ of $SO(4)$ which contain $H$. First we classify all the subgroups of $SU(2)\times SU(2)$ which contain  $\tilde{V}\times\tilde{V}$, then their images in $SO(4)$ are the groups $G$ (cf.\ proposition \ref{list} and section 1.4). They operate in a natural way on $\cix$, the ring of polynomials in four variables with complex coefficients. In section 3 we give generators for the spaces  $\cix^G_j$ of homogeneous $G$-invariant polynomials of degree $j$ whenever this dimension is two. Since the groups $G$ contain $H$, we have invariant polynomials only in even degree. When the dimension is two the generators are the multiple quadric \mbox{$q^{j/2}=(x_0^2+x_1^2+x_2^2+x_3^2)^{j/2}$} (trivial invariant) and another polynomial of degree $j$ which we denote by $f$. We describe  then the pencils
\begin{eqnarray*}
f+\lambda q^{j/2}=0,~~~~~\lambda\in\piu,
\end{eqnarray*}
 of symmetric surfaces in the three dimensional complex projective space $\pitr$ (section 4). In particular we find the singular surfaces contained in it.
In \cite{sa} we considered the case of $G=TT$, $OO$, $II$ which are the images in $SO(4)$ of the direct products $\tilde{T}\times\tilde{T}$, $\tilde{O}\times\tilde{O}$, $\tilde{I}\times\tilde{I}$ where $\tilde{T}$ denotes the binary tetrahedral group, $\tilde{O}$ the binary octahedral group, $\tilde{I}$ the binary icosahedral group in $SU(2)$. We denoted the groups there by $G_6$, $G_8$ and $G_{12}$  and we called them bi-polyhedral groups. We found pencils containing surfaces with many nodes (=ordinary double points). In particular the degree twelve $II$-invariant pencil contains a surface with $600$ nodes which improves the previous lower bound for the maximal number of nodes of a surface of degree twelve in $\pitr$ (cf. \cite{chmutov}). Here we describe the other $G$-invariant pencils and show that we have two more pencils which contain surfaces with nodes (the others do not contain surfaces with isolated singularities at all). We list the groups $G$ and the degrees $j$ below, as well as the number of nodes on the singular surfaces. In each pencil we have four of these singular surfaces and the nodes there form just one $G$-orbit. For the convenience of the reader we recall the results about the $TT$-, $OO$-, and $II$-invariant pencils too.
 
\renewcommand{\arraystretch}{1.3}

\begin{center}

\begin{tabular}{l|*{2}{l|}*{4}{l}}
$G$& {\rm order}& $j$ &\multicolumn{4}{c}{nodes}\\
\hline 
$(OO)'$&192&4&4&12&16&8\\
$TT$&288&6&12&48&48&12\\
$OO$&1152&8&24&72&144&96\\
$IO$&2880&12&240&360&240&120\\
$II$&7200&12&300&600&360&60
\end{tabular}

\end{center}
The group $IO$ is the image in $SO(4)$ of the direct product of the binary icosahedral group with the binary octahedral group in $SU(2)$$\times$$SU(2)$ and $(OO)'$ is a subgroup of $OO$ which we describe in section 1.2 and  1.4. This is an index two subgroup of the reflection group [3,3,4] (cf.\ \cite{cox}, p.\ 226) and has the same invariant polynomials. In all the cases but one ($G=IO$), the singular surfaces contain just isolated singularities (i.e.\ the nodes). The surfaces in the $IO$-invariant pencil contain two double lines in the base locus. \\
The [3,3,4]-invariant polynomials of degree two and four were already known by Coxeter in \cite{cox0}. Here we show that in the pencil of degree four we have a surface with the maximal number possible of nodes (=16) which is a so called Kummer surface. Finally in section 6 we give a computer picture of the $IO$-invariant surface of degree 12 with 360 nodes.\\
I thank Prof.\ Wolf Barth at the University of Erlangen for many helpful comments and discussions. 




\section{Symmetry groups}
Denote by $H\subseteq SO(4)$ the Heisenberg group (with 32 elements). We want to collect in a systematic way all the finite subgroups of $SO(4)$ containing $H$, and their polynomial invariants of low degree. These are invariants of the Heisenberg group with extra symmetries.
\subsection{Ternary groups}
We specify the following matrices in $SO(3)$
\begin{eqnarray*}
A_1:=\left(\begin{array} {ccc}
1&0&0\\
0&-1&0\\
0&0&-1
\end{array}\right ), A_2:=\left(\begin{array} {ccc}
-1&0&0\\
0&1&0\\
0&0&-1
\end{array}\right ),
R_n:=\left(\begin{array} {ccc}
1&0&0\\
0&a&-b\\
0&b&a
\end{array}\right ),\\
S:=\left(\begin{array} {ccc}
0&-1&0\\
0&0&-1\\
1&0&0
\end{array}\right ),U:=\frac{1}{2}\left(\begin{array} {ccc}
\tau-1&-\tau&1\\
\tau&1&\tau-1\\
-1&\tau-1&\tau
\end{array}\right ), 
\end{eqnarray*}
where $\tau:=\frac{1}{2}(1+\sqrt{5})=2\cdot$cos$(\frac{\pi}{5})$, $a:= $cos$\frac{2\cdot\pi}{n}$, $b:= $sin$\frac{2\cdot\pi}{n}$. These matrices generate the following subgroups of $SO(3)$
\begin{center}

\begin{eqnarray*}
\begin{array}{lllll}
&{\rm generators}&{\rm order}&{\rm group}& {\rm name}\\
V&A_1,A_2&4&\mathbb{Z}_2\times\mathbb{Z}_2&{\rm Klein~four}\\
D_n&A_2,R_n&2n&D_n&{\rm dihedral}\\
T&A_1,S&12&Alt(4)&{\rm tetrahedral}\\
O&A_1,R_4,S&24&Sym(4)&{\rm octahedral}\\
I&A_1,S,U&60&Alt(5)&{\rm icosahedral}
\end{array}
\end{eqnarray*}
\end{center}
Where $Alt(4)$ and $Alt(5)$ denote the group of even permutations of four  and five elements and $Sym(4)$ denotes the permutation group of four elements.
Whenever $n\in \mathbb{N}$, $n\not=0$, is {\it even}, $V$ is contained in each of the above groups. By the classification of the finite subgroups of $SO(3)$, these are all such subgroups. \\
The groups $T,O$ and $I$ are the rotation groups of tetrahedron, octahedron and icosahedron. An identification with the permutation groups is given in \cite{cox}, p.\ 46-50 as well as an identification of $D_n$ with the symmetry group of a regular polygon with $n$ vertices, p.\ 46. Sometimes it is useful for the computation to identify the matrices above with the cycles of the permutation groups. Indeed the identification of $T$ and $O$ with subgroups of the permutation  group $Sym(4)$ is obtained by letting them act on the four space diagonals of the unit cube. Let these lines and vectors generating them be
\begin{eqnarray*}
d_1:(1,1,1),~d_2:(-1,1,1),~d_3:(1,-1,1),~d_4:(1,1,-1)\\
\end{eqnarray*}
The matrices in $O$ permute these lines by
\begin{eqnarray*}
A_1:(12)(34),~A_2:(13)(24),~R_4:(1423),~S:(123).
\end{eqnarray*}
Using this correspondence with permutation groups, it is easy to write down their conjugacy classes. We write the conjugacy classes of the dihedral group $D_{n}$, $n=2l$, $l\geq 2$ too. In the next table we characterize a conjugacy class by one of its elements. Under each representative we write the number of elements in the conjugacy class.\\
\begin{center}
\begin{tabular}{l|lllll}
group &\multicolumn{5}{c}{repr. of a conj. class}\\
&\multicolumn{5}{c}{and number of elements}\\
\hline
$V$&1&$A_1$&$A_2$&$A_3$&\\
&1&1&1&1\\
\hline
$D_{n}$&1&$R_n^k$&$R_n^l$&$A_2$&$A_2R_n$\\
&1&2&1&l&l\\
\hline
T&1&$A_1$&$S$&$S^2$&\\
&1&3&4&4\\
\hline
$O$&1&$A_1$&$R_4A_2$&$R_4$&$S$\\
&1&3&6&6&8
\end{tabular}
\end{center}
where $k=1,\ldots,l-1$.\\
The symmetries of $T$ obviously leave invariant the icosahedron (cf.\ \cite{cox}, p.\ 52) with vertices 
\begin{eqnarray*}
(\pm 1,\pm \tau,0),~(0,\pm 1,\pm \tau),~(\pm \tau, 0,\pm 1)
\end{eqnarray*}
The matrix $U$ permutes these vertices as
\begin{eqnarray*}
\begin{array}{c}
\pm(0,1,\tau)\mapsto\pm(0,1,\tau),\\
\pm(0,1,-\tau)\mapsto\pm(-\tau,0,1)\mapsto\pm(-1,-\tau,0)\mapsto\pm(1,-\tau,0)\mapsto\pm (\tau,0,-1)\mapsto\pm(0,1,-\tau).
\end{array}
\end{eqnarray*}
So together with the group $T$ the symmetry $U$ generates a group of order at least $12\cdot 5=60$ , contained in the symmetry group of the icosahedron specified. Therefore it coincides with the icosahedral group $I\cong A_5$. The action of $I$ on the five cosets of $I/T$ defines the map $I\mapsto A_5$. Using
\begin{eqnarray*}
\begin{array}{c}
A_1\cdot U=U^4\cdot A_1,~A_1\cdot U^2=U^3\cdot A_1,\\
S\cdot U=U^3\cdot A_2,~S\cdot U^2=U\cdot A_2S,~S\cdot U^3=U^2\cdot S^2,~S\cdot U^4=U^4\cdot S^2 A_2,\\
\end{array}
\end{eqnarray*}
one finds that under this map 
\begin{eqnarray*}
A_1\mapsto(14)(23),~S\mapsto(132),~U\mapsto(12345).
\end{eqnarray*}
Using this correspondence, one enumerates the conjugacy classes in $I$\\
\begin{center}
\begin{tabular}{l|lllll}
group &\multicolumn{5}{c}{repr. of a conj. class}\\
&\multicolumn{5}{c}{and number of elements}\\
\hline
$I$&1&$A_1$&$S$&$U$&$U^2$\\
&1&15&20&12&12
\end{tabular}
\end{center}
\subsection{Subgroups of products of ternary groups}
Here we classify {\it subdirect products} $G\subseteq G_1\times G_2$, where $G_1$ and $G_2$ are finite ternary groups $V,D_n$ ($n$ even), $T,O$ or $I$. We assume that $G$ contains the subgroup $V\times V$$\subseteq G_1\times G_2$. We are interested in these subgroups only up to interchanging the factors $G_1$ and $G_2$. So we assume that we are in one of the following cases
\begin{itemize}
\item $|G_1|\geq|G_2|$ and $G_1,G_2\not=D_n$,
\item $G_1\not=D_n$, $G_2=D_n$, 
\item $G_1=D_n$ and $G_2=D_m$.\\
\end{itemize}
Additionally, passing to smaller subgroups $G_i'$$\subseteq$$G_i$ if necessary, we may assume that both projections 
\begin{eqnarray*}
p_1:G\rightarrow G_1',~p_2:G\rightarrow G_2'
\end{eqnarray*}
are surjective, so that $G\subseteq G_1'\times G_2'$. Finally, we do not distinguish between groups conjugate in $SO(3)\times SO(3)$. 
In the table below we 
assume $n\not=m$, $n=2l$, $m=2l'$ and let $s:=${\rm lcm}$(n,m)$.\\
\begin{prop}\label{list}
The following list is a complete list of subgroups $G\subseteq G_1\times G_2$ under the assumptions above:
\pagebreak
\begin{eqnarray*}
\begin{array}{l|l|lll|l|l|lll}
G_1&G_2&G&|G|/16&G/V\times V&G_1&G_2&G&|G|/16&G/V\times V\\
\hline
V&V&V\times V&1&1&V&D_n&V\times D_n&n/2&\zi_l\times 1\\
T&V&T\times V&3&\zi_3\times 1&T&D_n&T\times D_n&3n/2&\zi_3\times\zi_l\\
&T&T\times T&9&\zi_3\times\zi_3&& 3|n &(T\times D_n)'&n/2&\zi_l\\
&&(T\times T)'&3&\zi_3&O&D_n&O\times D_n&6n/2&D_3\times\zi_l\\
O&V&O\times V&6&D_3\times1&&4|n&(O\times D_n)'&3n/2&D_{3l/2}\\
&T&O\times T&18&D_3\times\zi_3&I&D_n& I\times D_n&15n/2&\\ 
&O&O\times O&36&D_3\times D_3&D_n&D_n&D_n\times D_n&n^2/4&\zi_l\times\zi_l\\
&O&O\times O&36&D_3\times D_3&&&(D_n\times D_n)'&n/2&\zi_l\\
&&(O\times O)'&6&D_3&D_n&D_m&D_n\times D_m&nm/4&\zi_l\times\zi_{l'}\\
&&(O\times O)''&18&\zi_3\times\zi_3\times\zi_2&&&(D_n\times D_m)'&s/2&\zi_s\\
I&V&I\times V&15&&&&&&\\
&T&I\times T&45&&&&&&\\
&O&I\times O&90&&&&&&\\
&I&I\times I&225&&&&&&\\
\end{array}
\end{eqnarray*}
\end{prop}
\bprf
We discuss the cases one-by-one. But before that, we observe that the kernels $K_1\subseteq G_1\times 1$ and $K_2\subseteq 1\times G_2$ of both projections $p_i:G\rightarrow G_i$ are normal subgroups. This follows by conjugating component-wise from the surjectivity of both projections.\\
The groups $G_1\times V$ do not have proper subgroups containing $V\times V$ and mapping surjectively onto $G_1$, (for $G_1=D_n$ too). So we do not need to consider the cases $G_2=V$.\\
First consider the case of $G_1=I$. Since $I$ is simple, the kernel $K_1\subseteq I\times 1$ either coincides with $I\times 1$, or is trivial. The latter cannot be the case, because this kernel contains $V\times 1$. The only possibilities are the product cases $I\times V$, $I\times T$, $I\times O$, $I\times I$ and $I\times D_n$.\\
Let now $Q:=G/V\times V$ and $Q_1:=G_1/V$, $Q_2:=G_2/V$. We consider $Q$ a proper subdirect product of $Q_1\times Q_2$.\\
$T,T$: $Q\subseteq \zi_3\times\zi_3$ mapping surjectively onto both factors is either the diagonal or the anti-diagonal.
The two corresponding subgroups are not conjugate in $T\times T$, but in $T\times O$ they are. The inverse image of the diagonal $\zi_3$$\subseteq\zi_3\times\zi_3$ in $T\times T$ is $(T\times T)'$.\\
$O,T$: $Q\subseteq D_3\times \zi_3$ mapping surjectively onto both factors would have order six and be isomorphic with $D_3$ under $p_1$. But there is no epimorphism of $D_3$ onto $\zi_3$. Such a group does not exist.\\
$O,O$: $Q\subseteq D_3\times D_3$ must have order six, twelve or 18. If it has order six, it is a graph of an isomorphism between both factors $D_3$. Then it is conjugate to the diagonal, and this leads to the subgroup $(O\times O)'$$\subseteq$$O\times O$. The case $|Q|$=$12$ cannot occur, because then the kernel $K_1$$\subseteq$$D_3\times 1$ would have order two, and could not be normal. If $Q$ has order $18$ both the kernels $K_1$ and $K_2$ have order three, and coincide with the unique proper normal subgroup of $D_3$.  This implies that $G$ contains $T\times T $$\subseteq$$O\times O$ and is the inverse image of a subgroup $\zi_2$$\subseteq$$\zi_2\times\zi_2$$=O\times O/T\times T$. By surjectivity of projections this can only be the diagonal. Its inverse image is $(O\times O)''$. \\
$T,D_n$: if $Q\subseteq \zi_3\times \zi_l$ maps surjectively onto both the components then three divides $l$. Let its inverse image be $(T\times D_n)'$.\\
$O,D_n$: if $Q\subseteq D_3\times \zi_l$ maps surjectively onto both factors then two divides $l$. We denote the inverse image by $(O\times D_n)'$.\\
$D_n,D_n$: $Q\subseteq \zi_l\times\zi_l$ mapping surjectively onto $\zi_l$ is conjugate to the diagonal. We call its inverse image $(D_n\times D_n)'$.\\
$D_n,D_m$:  $Q\subseteq\zi_l\times\zi_{l'}$ mapping surjectively onto $\zi_l$ and $\zi_{l'}$ is generated by an element of order $s/2$.
We call its inverse image $(D_n\times D_m)'$.\\
\eprf 
\subsection{Binary groups}
We consider the double cover of \cite{tits} p.\ 77--78. Let $\tilde{G}$ denotes the pre-image in $SU(2)$ of $G$$\subseteq$$SO(3)$. We specify the following matrices, $\tilde{M}$, which are in the pre-image of $M$$\in$$SO(3)$:
 \begin{eqnarray*}
\begin{array}{lll}
\tilde{A_1}:=\left(\begin{array} {cc}
i&0\\
0&-i
\end{array}\right ), &\tilde{A_2}:=\left(\begin{array} {cc}
0&1\\
-1&0
\end{array}\right ),& \tilde{A_3}:=\left(\begin{array} {cc}
0&i\\
i&0
\end{array}\right ),\\
\tilde{S}:=\frac{1}{2}\left(\begin{array} {cc}
1+i&-1+i\\
1+i&1-i
\end{array}\right ),& \tilde{U}:=\frac{1}{2}\left(\begin{array} {cc}
\tau&\tau-1+i\\
1-\tau+i&\tau
\end{array}\right ),& \tilde{R_n}:=\left(\begin{array} {cc}
e^{\frac{i\pi}{n}}&0\\
0&e^{-\frac{i\pi}{n}}
\end{array}\right ).
\end{array}
\end{eqnarray*}
since $\tilde{M}^{\mbox{ord}(M)}$$=-1$, they have order $2\cdot$ord$(M)$. By an argumentation as in \cite{sa} section 2, we can write the conjugacy classes in the binary groups:
\pagebreak
\begin{center}
\begin{tabular}{c|ccccccccc}
group &\multicolumn{9}{c}{repr. of a conj. class and number of elements}\\
\hline
$\tilde{V}$&1&-1&$\tilde{A_1}$&$\tilde{A_2}$&$\tilde{A_3}$&&&&\\
&1&1&2&2&2&&&&\\
\hline
$\tilde{D_n}$&1&-1&$\tilde{A_2}$&$\tilde{A_2}\tilde{R_n}$&$\tilde{R_n}^l$&$\tilde{R_n}^{k}$&$-\tilde{R_n}^k$&&\\
&1&1&2l&2l&2&2&2&&\\
\hline
$\tilde{T}$&1&-1&$\tilde{A_1}$&$\tilde{S}$&$-\tilde{S}$&$\tilde{S}^2$&$-\tilde{S}^2$&&\\
&1&1&6&4&4&4&4&&\\
\hline
$\tilde{O}$&1&-1&$\tilde{A_1}$&$\tilde{R_4}\tilde{A_2}$&$\tilde{R_4}$&$-\tilde{R_4}$&$\tilde{S}$&$-\tilde{S}$\\
&1&1&6&12&6&6&8&8&\\
\hline
$\tilde{I}$&1&-1&$\tilde{A_1}$&$\tilde{S}$&$-\tilde{S}$&$\tilde{U}$&$-\tilde{U}$&$\tilde{U}^2$&$-\tilde{U}^2$\\
&1&1&30&20&20&12&12&12&12
\end{tabular}
\end{center}
where $k=1,\ldots,l-1$.
\subsection{Quaternary groups}
Here we consider the images of the finite groups $\tilde{G_1}\times\tilde{G_2}\subseteq SU(2)\times SU(2)$ under the double covering map
\begin{eqnarray*}
\begin{array}{cc}
SU(2)\times SU(2)\rightarrow SO(4),&(q,q'):p\mapsto q\cdot p\cdot q^{-1}\\
\end{array}
\end{eqnarray*}
(cf.\ \cite{tits} p.\ 77--78), we abbreviate there $G_1G_2$. Since the corresponding subgroups of $SO(3)\times SO(3)$ contain the group $V\times V$, these contain  the Heisenberg group $VV$$\subseteq$$SO(4)$, and by proposition \ref{list} these are all such subgroups. We specify now the matrices:\\
\begin{eqnarray*}
\begin{array}{ll}
(A_1,1):=\left(\begin{array} {cccc}
0&-1&0&0\\
1&0&0&0\\
0&0&0&-1\\
0&0&1&0
\end{array}\right ),& (1,A_1):=\left(\begin{array} {cccc}
0&1&0&0\\
-1&0&0&0\\
0&0&0&-1\\
0&0&1&0
\end{array}\right ),\\
(A_2,1):=\left(\begin{array} {cccc}
0&0&-1&0\\
0&0&0&1\\
1&0&0&0\\
0&-1&0&0
\end{array}\right ),& (1,A_2):=\left(\begin{array} {cccc}
0&0&1&0\\
0&0&0&1\\
-1&0&0&0\\
0&-1&0&0
\end{array}\right ),\\
\end{array}
\end{eqnarray*}
\begin{eqnarray*}
\begin{array}{ll}
(A_3,1):=\left(\begin{array} {cccc}
0&0&0&-1\\
0&0&-1&0\\
0&1&0&0\\
1&0&0&0
\end{array}\right ), &(1,A_3):=\left(\begin{array} {cccc}
0&0&0&1\\
0&0&-1&0\\
0&1&0&0\\
-1&0&0&0
\end{array}\right ),\\
(R_n,1):=\left(\begin{array} {cccc}
\alpha&-\beta&0&0\\
\beta&\alpha&0&0\\
0&0&\alpha&-\beta\\
0&0&\beta&\alpha
\end{array}\right ), &(1,R_n):=\left(\begin{array} {cccc}
\alpha&\beta&0&0\\
-\beta&\alpha&0&0\\
0&0&\alpha&-\beta\\
0&0&\beta&\alpha
\end{array}\right ),\\
(S,1):=\frac{1}{2}\left(\begin{array} {cccc}
1&-1&1&-1\\
1&1&-1&-1\\
-1&1&1&-1\\
1&1&1&1
\end{array}\right ), &(1,S):=\frac{1}{2}\left(\begin{array} {cccc}
1&1&-1&1\\
-1&1&-1&-1\\
1&1&1&-1\\
-1&1&1&1
\end{array}\right ),\\
(U,1):=\frac{1}{2}\left(\begin{array} {cccc}
\tau&0&1-\tau&-1\\
0&\tau&-1&\tau-1\\
\tau-1&1&\tau&0\\
1&1-\tau&0&\tau
\end{array}\right ), &(1,U):=\frac{1}{2}\left(\begin{array} {cccc}
\tau&0&\tau-1&1\\
0&\tau&-1&\tau-1\\
1-\tau&1&\tau&0\\
-1&1-\tau&0&\tau
\end{array}\right ),
\end{array}
\end{eqnarray*}
where $\alpha:=$ cos$\frac{\pi}{n}$, $\beta:=$ sin$\frac{\pi}{n}$.\\
We can write the conjugacy classes of the groups $G_1G_2$$\subseteq$$SO(4)$ in $SO(4)$ and their number of elements. These are the images of the conjugacy classes of $\tilde{G_1}\times\tilde{G_2}$ in $G_1G_2$$\subseteq$$SO(4)$ under the double covering map. Since the fact (1.1) of \cite{sa} still holds, the matrices $(g_1,g_2)$$\in SO(4)$ with the same eigenvalues are conjugate. This simplifies the computations considerably. In this section and in the next one we omit the groups $D_nD_m$ and $(G_1G_2)'$ with $G_2=D_n$. We return to those groups later.

In the tables we use the following conventions:
\begin{itemize}
\item we omit the conjugacy classes of $+1, -1$ (these contain one element each)
\item whenever the conjugacy classes $(g_1,g_2)$ and its $s$th power $(g_1^s,g_2^s)$ are distinct we write them just one time since they have the same number of elements.
\end{itemize}
 
\begin{eqnarray*}
\begin{array}{c|c|ccccccccc}
G_1G_2&{\rm order}&A_2,1&A_2,A_2&R_4,1&R_4,A_2&R_4,R_4&S,1&S,A_2&S,R_4&S,S\\
\hline
VV&32&12&18&0&0&0&0&0&0&0\\
TV&96&12&18&0&0&0&8&48&0&0\\
TT&288&12&18&0&0&0&16&96&0&64\\
OV&192&24&54&6&36&0&8&48&0&0\\
OT&576&24&54&6&36&0&16&192&48&64\\
OO&1152&36&162&12&216&36&16&288&96&64\\
IV&480&36&90&0&0&0&20&120&0&0\\
IT&1440&36&90&0&0&0&28&360&0&160\\
IO&2880&48&270&6&180&0&28&600&120&160\\
II&7200&60&450&0&0&0&40&1200&0&400\\
\end{array}
\end{eqnarray*}
\begin{eqnarray*}
\begin{array}{c|c|cccccc}
G_1G_2&{\rm order}&A_2,1&A_2,A_2&R_4,1&R_4,A_2&S,1&S,A_2\\
\hline
VD_n&16n&8+4l&6(2l+1)&0&0&0&0\\
TD_n&48n&8+4l&6(2l+1)&0&0&8&16(2l+1)\\
OD_n&96n&20+4l&18(2l+1)&6&12(2l+1)&8&16(2l+1)\\
ID_n&240n&32+4l&30(2l+1)&0&0&20&40(2l+1)\\
\end{array}
\end{eqnarray*}
\begin{eqnarray*}
\begin{array}{c|ccccccc}
G_1G_2&U,1&U,A_2&U,R_4&U,S&U,U&U,U^2&U,R_n^k\\
\hline
IV&12&72&0&0&0&0&0\\
IT&12&72&0&96&0&0&0\\
IO&12&216&72&96&0&0&0\\
II&24&720&0&480&144&144&0\\
ID_n&12&24(2l+1)&0&0&0&0&24\\
\end{array}
\end{eqnarray*}

\begin{eqnarray*}
\begin{array}{c|cccc}
G_1G_2&A_2,R_n^k& 1,R_n^k&S,R_n^k&R_4,R_n^k\\
\hline
VD_n&12&2&0&0\\
TD_n&12&2&16&0\\
OD_n&36&2&16&12\\
ID_n&60&2&40&0\\
\end{array}
\end{eqnarray*}
where $k=1,\ldots,l-1$. 
For the groups which are not products it is a little more complicated to write down the sizes of their conjugacy classes. But using the description from proposition \ref{list} one finds\\
\begin{eqnarray*}
\begin{array}{c|c|ccccccc}
{\rm group}&{\rm order}&A_2,1&A_2,A_2&R_4,A_2&R_4,R_4&S,1&S,A_2&S,S\\
\hline
(TT)'&96&12&18&0&0&0&0&32\\
(OO)'&192&12&42&48&12&0&0&32\\
(OO)''&576&12&90&144&36&16&96&64\\
\end{array}
\end{eqnarray*}


\section{Poincar\'e series}
In this section we want to find the dimension of the spaces of homogeneous invariant polynomials of a given degree. We consider the Poincar\'e series 
\begin{eqnarray*}
p(t):=\sum_{j=0}^{\infty}{\rm dim}\cix^{G}_j \cdot t^j
\end{eqnarray*}
 where $G$ is a group as in section  1.4. By a theorem of Molien (\cite{ben} p.\ 21) and an easy computation as in \cite{sa}  (2.1), it  can be written as
\begin{eqnarray*}
p(t)=\frac{1}{|G|}\sum\frac{n_g}{{\rm det}(g-1\cdot t)}
\end{eqnarray*}
where the sum runs over all the conjugacy classes of $G$ and $n_g$ denote their number of elements. At the denominator we have the characteristic polynomials. Using the numbers of conjugacy classes (under $SO(4)$) given in section 1.4 and computing their characteristic polynomials, the power series package of MAPLE produces the following table of dimensions $m_d$ of invariant polynomials in degree $d$ of the groups $G$$\subseteq$$G_1G_2$. Observe that since $G$ contains the Heisenberg group we do not have invariant polynomials of odd degree. First we consider the case of $G_i$$\not=D_n$:
\begin{eqnarray*}
\begin{array}{c|cccccc}
{\rm group}&m_2&m_4&m_6&m_8&m_{10}&m_{12}\\
\hline
VV&1&5&6&15&19&35\\
TV&1&1&{\bf 2}&5&5&13\\
TT&1&1&{\bf 2}&3&3&7\\
(TT)'&1&3&4&7&9&15\\
OV&1&1&1&4&4&8\\
OT&1&1&1&{\bf 2}&{\bf 2}&4\\
OO&1&1&1&{\bf 2}&{\bf 2}&3\\
(OO)'&1&{\bf 2}&3&5&6&9\\
(OO)''&1&1&{\bf 2}&3&3&5\\
IV&1&1&1&1&1&5\\
IT&1&1&1&1&1&3\\
IO&1&1&1&1&1&{\bf 2}\\
II&1&1&1&1&1&{\bf 2}\\
\end{array}
\end{eqnarray*}
Whenever $G_1\not=D_n$ and  $G_2=D_n$, we can write the finite sums $p(t)$ for each $n$. Here we compute the first coefficients of the Poincar\'e series for  $n=4,6,8$. 
\begin{eqnarray*}
\begin{array}{c|cccccc}
{\rm group}&m_2&m_4&m_6&m_8&m_{10}&m_{12}\\
\hline
VD_4&1&3&3&9&11&19\\
VD_6&1&3&3&6&6&14\\
VD_8&1&3&3&6&6&10\\
\hline
TD_4&1&1&1&3&3&7\\
TD_6&1&1&1&{\bf 2}&{\bf 2}&6\\
TD_8&1&1&1&{\bf 2}&{\bf 2}&4\\
\hline
OD_4&1&1&1&3&3&5\\
OD_6&1&1&1&{\bf 2}&{\bf 2}&4\\
OD_8&1&1&1&{\bf 2}&{\bf 2}&3\\
\hline
ID_4&1&1&1&1&1&3\\
ID_6&1&1&1&1&1&3\\
ID_8&1&1&1&1&1&{\bf 2}
\end{array}
\end{eqnarray*}
Of course, whenever $m_d=1$, the space of invariant polynomials is spanned by the $d/2$-th power of the invariant quadric ({\it trivial invariant})
\begin{eqnarray*}
q:=x_0^2+x_1^2+x_2^2+x_3^2.
\end{eqnarray*}

\input{artjul01_1}

%% file: artjul01_1.tex
\section{Invariants} \label{invariants}
In this section we want to compute a system of generators for the spaces $\cix_j^G$ whenever this space has dimension two and $G$$\subseteq$SO$(4)$ is a finite subgroup containing the Heisenberg group. We distinguish two cases.\\ 
{\bf (3.1)}{\it First case}. Assume that $G\subseteq G_1G_2$ with $G_i$$\not=D_n$, $i=1,2$. We do some remarks on the groups $G$ which simplify the computations of the invariant polynomials.
\begin{itemize}
\item the group $TV$$\subseteq$$TT$.
\item The group $TT$ is contained in $(OO)''$. In fact the generators modulo $VV$ are 
\begin{eqnarray*}
\begin{array}{l|ll}
{\rm group} &TT&(OO)''\\
\hline
{\rm generators}&(1,S)&(1,S),(S,1)\\
&(S,1)&(R_4A_2,R_4A_2)
\end{array}
\end{eqnarray*}
with 
\begin{eqnarray*}
(R_4A_2,R_4A_2)=\left(\begin{array} {cccc}
1&0&0&0\\
0&-1&0&0\\
0&0&0&1\\
0&0&1&0
\end{array}\right )
\end{eqnarray*} 
Similarly to \cite{sa} section  3., $(OO)''$ is a subgroup of index two in the reflection group of the $\{3,4,3\}$-cell (cf.\ \cite{cox} p.\ 149 for the definition of this polytope), if we add the generator 
\begin{eqnarray*}
C=\left( \begin{array} {cccc}
1 & 0& 0 & 0 \\
0 & -1& 0 & 0 \\
0 & 0& -1& 0 \\
0 & 0& 0 & -1 
\end{array} \right)
\end{eqnarray*}
 we get the whole reflection group. 
\item Finally since $OT$$\subseteq$$OO$, they have the same two invariant polynomials in degree eight. Those of degree ten are obtained just by multiplication with the quadric $q$.  
\end{itemize}

By these remarks and considering $G$ as in the assumption, it follows that we have to compute the generators of six two-dimensional spaces $\cix_j^G$ for the following pairs $(G,j)$:
\begin{eqnarray*}
\begin{array}{lllll}
((OO)',4)&(TT,6)&(OO,8)&(IO,12)&(II,12).\\
\end{array}
\end{eqnarray*}
The generators of the $TT$-, $OO$- and $II$-invariant spaces are given in \cite{sa} as well as a description of the corresponding pencils of surfaces in $\pitr$ (base locus and singular surfaces). Here we examine the remaining  cases and in the next section we describe the corresponding pencils of surfaces in $\pitr$. The basic idea to find generators of the invariant spaces is the same as in \cite{sa}. \\
{\it $(OO)'$-invariants}.
We start with the space of Heisenberg invariant quartics. It has dimension five, being spanned by 
\begin{eqnarray*}
\begin{array}{lll}
f_0:=x_0^4+x_1^4+x_2^4+x_3^4,&&\\
f_1:=2(x_0^2x_1^2+x_2^2x_3^2),&f_2:=2(x_0^2x_2^2+x_1^2x_3^2),&f_3:=2(x_0^2x_3^2+x_1^2x_2^2),\\
f_4:=4x_0x_1x_2x_3.&&\\
\end{array}
\end{eqnarray*}
In terms of these invariants
\begin{eqnarray*}
q^2=f_0+f_1+f_2+f_3.
\end{eqnarray*}
Modulo $VV$, the group $(OO)'$ is generated by $(R_4A_2,R_4A_2)$ and $(S,S)$.
Tracing the action of these generators on $f_0,\ldots, f_4$ one finds the invariants $q^2$ and $f_0$.\\

{\it $IO$-invariants}. 
The generators $(U,1)$, $(1, S)$ and $(1, R_4)$ of $IO$ operate on the space $\cix^{VV}_{12}$ which is 35-dimensional. A computation with MAPLE shows that it  contains the non trivial $IO$-invariant polynomial
\begin{eqnarray*}
\begin{array}{lll}
S_{IO}&:=&-\sum_{i} f_i^3+11\sum_{i,j} f_i^2f_j+(f_0^2-14f_4^2)\sum_{i}f_i+30f_0f_4^2-2f_0\sum_{i,j}f_if_j\\
&&-30f_1f_2f_3+3\sqrt{5}f_4(2f_0\sum_{i}f_i-\sum_{i,j}f_if_j-f_0^2-\sum_{i}f_i^2+4f_4^2)\\
&&+6\sqrt{5}F_a
\end{array}
\end{eqnarray*}
where the sums run over all the indices $i,j=1,2,3$, $i\not=j$, and 
\begin{eqnarray*}
F_a:=f_1^2f_2+f_3^2f_1+f_2^2f_3-f_1f_2^2-f_3f_1^2-f_2f_3^2
\end{eqnarray*}
is the anti-symmetric part of $S_{IO}$.\\
{\bf(3.2)}{\it Second case}. Assume that $G$$\subseteq$$G_1G_2$ with $G_1=D_n$ or $G_2=D_m$, $n,m\geq 4$, even. With the help of the table in section 1.4 we discuss the following cases:
\begin{itemize}
\item $D_nD_m$: observe that $VD_m$ is contained in $D_nD_m$. A direct computation shows that the generator $(1,R_n)$ let invariant a three-dimensional family of degree four $VV$-invariant polynomials. Generators in this case are $f_2+f_3$, $q^2$, $f_2-f_4$. Now the matrix $(R_n,1)$ of $D_nD_m$ operates on these three-dimensional space leaving invariant $q^2$ and $f_2+f_3$.
\item since $(D_nD_m)'$$\subseteq$$D_nD_m$ these groups have already an at least two-dimensional family of invariant polynomials of degree four, which is generated by $q^2$ and $f_2+f_3$.
\item $TD_n$, $n\geq 6$: we have a five-dimensional family of $TV$-invariant in degree eight. By the action of the generator $(1,R_n)$, $n\geq 6$ we have a two-dimensional family of invariants. Put 
\begin{eqnarray*}
\begin{array}{l}
K_i:=x_0^2+x_1^2-x_2^2-x_3^2+(-1)^i2x_0(x_2+x_3)+2x_1(x_2-x_3),~i=0,3\\
K_i:=x_0^2+x_1^2-x_2^2-x_3^2+(-1)^i2x_0(x_2+x_3)-2x_1(x_2-x_3),~i=1,2
\end{array}
\end{eqnarray*}
then generators are 
\begin{eqnarray*}
q^4 ~\mbox{and}~ P_8:=\prod_{i=0}^{3}K_i.
\end{eqnarray*}
\item Suppose three divides $n$. Since the groups  $(TD_n)'$ are contained in $TD_n$, they have already at least a two-dimensional family of invariant polynomials in degree eight. We consider now the action of the extra generators $(S,R_n)$ of $(TD_n)'$ on the space of degree four, resp. six $VV$-invariant polynomials. A direct computation with MAPLE shows that we have no invariants other then the quadric $q^2$, resp. $q^3$. 
\item $OD_n$, $n\geq 6$: we have $TD_n$$\subseteq$$OD_n$ and the extra generator $(R_4,1)$ leaves the previous polynomials invariant.
\item Suppose now four divides $n$. Since the groups $(OD_n)'$ are contained in $OD_n$, they have at least a two-dimensional family of invariant polynomials in degree eight. Modulo $VV$ the groups $(OD_n)'$ have generators $(S,1)$ and $(R_4,R_n)$, hence they contains the group $TV$.  This has no non-trivial invariant polynomials of degree four and has a two-dimensional family of invariant polynomials in degree six, generated by $q^3$ and $S_6(x)$ (cf. \cite{sa} p.\ 437). For $n\geq 6$, $S_6(x)$ is not $(OD_n)'$-invariant, for $n=4$ it is. In any case we do not get new invariants.  
\item $ID_n$, $n\geq 8$: we have a five-dimensional family of $IV$-invariant polynomials in degree twelve. The action of $(1,R_n)$, $n\geq 8$ on this space produces a two-dimensional family of invariant polynomials. Put 
\begin{eqnarray*}
\begin{array}{l}
K_0':=x_0^2+x_1^2-x_2^2-x_3^2+2(1-\tau)(x_0x_2-x_1x_3)\\
K_1':=x_0^2+x_1^2-x_2^2-x_3^2+2\tau(x_0x_3+x_1x_2)\\
K_2':=x_0^2+x_1^2-x_2^2-x_3^2-2(1-\tau)(x_0x_2-x_1x_3)\\
K_3':=x_0^2+x_1^2-x_2^2-x_3^2-2\tau(x_0x_3+x_1x_2)\\
K_4':=2(1-\tau)(x_0x_3+x_1x_2)+2(x_0x_2-x_1x_3)\\
K_5':=-2\tau(x_0x_3+x_1x_2)+2(x_0x_2-x_1x_3)
\end{array}
\end{eqnarray*}
then generators are 
\begin{eqnarray*}
q^6~\mbox{and}~P_{12}:=\prod_{i=0}^{5}K_i'.
\end{eqnarray*}
\end{itemize}
In conclusion, the dimension of $\cix^G_j$ is two, for $(G,j)$ equal to
 \begin{eqnarray*}
\begin{array}{lllll}
(D_nD_m,4)~ n,m\geq 4,&(OD_n,8)~n\geq 6,&(ID_n,12)~n\geq 8.\\
\end{array}
\end{eqnarray*}
{\bf Remark 3.3} Observe that the polynomials $q^2$ and $f_2+f_3$ are $D_nD_n$-invariant even if $n$ is an odd integer. As in the even case, generators of these groups are $(A_2,1)$, $(1,A_2)$ and $(R_n,1)$, $(1,R_n)$.

\section{Invariant pencils}
\subsection{The pencil of $(OO)'$-invariant quartics}
We take the generators $q^2$ and $f_0:=x_0^4+x_1^4+x_2^4+x_3^4$. The base locus is a double curve of degree eight. Let $C_8$ denote the curve $\{q=0,~f_0=0\}$. The pencil is invariant under the matrix $C$ too, hence by symmetry reasons the curve $C_8$ has bi-degree $(4,4)$ on $q$ (cf.\ \cite{sa}, section 5). Moreover we have the following
\begin{lemma} The curve $C_8$ is smooth and irreducible. \end{lemma}
{\it Proof.} The Jacobian matrix of $C_8$ has rank two in each point, hence the curve is smooth. If $C_8=C'\cup C''$, the curves $C'$, $C''$ would meet in some point on $q$ (observe that they cannot be lines of the same ruling), but this is impossible because $C_8$ is smooth. \eprf
In particular (6.1) of \cite{sa} still holds and by Bertini's theorem the general surface in the pencil is smooth. All the other surfaces but $q^2$ are irreducible and reduced and the singular ones have only isolated singularities. \\
{\it The symmetry group of $\{3,3,4\}$}. Consider the ``cross polytope'' $\beta_4$$=\{3,3,4\}$ in $\mathbb{R}^4$ with vertices the permutations of $(\pm 1, 0, 0,0)$ as in \cite{cox} p.\ 156 and edge $\sqrt{2}$. The generators of $(OO)'$ permutes these points, hence $(OO)'$ is contained in the symmetry group $[3,3,4]$ of $\{3,3,4\}$.  More precisely it is an index two subgroup, in fact the symmetry group of $\beta_4$ has order $2^4\cdot 4!=384$$=2\cdot 192$ and by adding the generator $C$ to $(OO)'$ we get the whole symmetry group $[3,3,4]$ (cf. \cite{cox} p.\ 226). In particular observe that they have the same invariant polynomials. The polytope $\beta_4$ has $N_0=8$, $N_1=24$, $N_2=32$, $N_3=16$ the reciprocal ``measure polytope'', $\gamma_4$$=\{4,3,3\}$ has $N_0=16$, $N_1=32$, $N_2=24$, $N_3=8$. Hence we get four $[3,4,3]$-orbits of points: the vertices and the middle points of the edges of the $\beta_4$ and the vertices and the middle points of the edges of the  reciprocal $\gamma_4$. These have coordinates the permutation of $(\pm 1,0,0,0)$, $(\pm 1,\pm 1, 0,0)$,resp. $(\pm 1,\pm 1,\pm 1, \pm 1)$, $(\pm 1, \pm 1, \pm 1, 0)$. As points of $\pitr$ these are singular on the surfaces $f_0+\lambda q^2$ for $\lambda=-1$, $-\frac{1}{2}$ resp. $-\frac{1}{4}$, $-\frac{1}{3}$ and a direct computation shows that they are all ordinary double points. We do it for $\lambda=-1$ and $(1:0:0:0)$. In the affine chart $\{x_0\not=0\}$ the equation becomes 
\begin{eqnarray*}
\begin{array}{lll}
0&=&1+x^4+y^4+z^4-(1+x^2+y^2+z^2)^2\\
&=&1-2x^2-2y^2-2z^2+{\rm terms~of~degree} \geq 4,
\end{array}
\end{eqnarray*}
hence the rank of the Hessian matrix at $(0,0,0)$ is three. This shows that $(1:0:0:0)$ is an ordinary double point and so are all the points in its orbit.\\  We collect the results on the singular surfaces in the following table. In the middle column we write just one point, but we mean all its permutations.

\begin{eqnarray*}
\begin{array}{r|c|l}
\lambda& {\rm nodes}& {\rm number}\backslash {\rm description}\\ 
\hline
-1&(1:0:0:0)&4\\
\hline
-\frac{1}{2}&(\pm 1:\pm 1:0:0)&12\\
\hline
-\frac{1}{3}&(\pm 1:\pm 1:\pm 1:0)&16, {\rm Kummer~surface}\\
\hline
-\frac{1}{4}&(\pm 1:\pm 1:\pm 1:\pm 1)&8\\
\hline
\hline
\infty&-&{\rm double~quadric}
\end{array}
\end{eqnarray*}

In the case of $\lambda=-\frac{1}{3}$, we get a surface with 16 nodes which is the maximal number possible for a surface of degree four. This is a Kummer surface. Observe that in this case too as in \cite{sa} the nodes are fix points under the action of some matrices in $[3,3,4]$, resp. in $(OO)'$ and they are contained on lines of fix points, see  (5.2) and  (6.3) of \cite{sa}. Moreover an estimation as in section 8 of \cite{sa} shows that whenever the lines of fix points do not meet the base locus, they contain exactly four nodes. It is natural to expect that these are all the singular surfaces in the $(OO)'$-invariant pencil (as in the case of the $TT$-, $OO$-, and $II$-invariant pencils). This is a direct consequence of the following\\
\begin{prop} 
1. The conjugacy classes in $(OO)'$ (under $(OO)'$) with eigenvalues $\pm 1$ are the following
\begin{eqnarray*}
\begin{array}{l|c|c|c|c|c}
{\it conj.~class}&(A_2,A_2)&(A_2,A_1)&(R_4A_2,R_4A_2)&(S,S)&(R_4,R_4)\\
\hline
{\it number~of~elements}&6&12&24&32&12\\
{\it number~of~fix~lines}&6&12&24&16&6\\
\end{array}
\end{eqnarray*}
2. The fix lines of the matrices in these conjugacy classes contain the maximal number possible of node.\\
\end{prop}
{\it Proof.} Choosing a fix line for each of the representative above and intersecting with the singular surfaces we find
\vspace*{0.3cm}

\begin{eqnarray*}
\begin{array}{c|c|c|c|c}
{\rm matrix~and~fix~line}&\multicolumn{2}{c|}{{\rm value~of} ~\lambda}&\multicolumn{2}{c}{{\rm int.~points}}\\
\hline
(A_2,A_2): &-1&-\frac{1}{2}&(1:0:0:0)&(\pm 1:0:1:0)\\
x_1=x_3=0&&&(0:0:1:0)&\\
(A_2,A_1):&-\frac{1}{2}&-\frac{1}{4}&(0:1:1:0)&(1:1:1:1)\\
x_0=x_3, x_1=x_2&&&(1:0:0:1)&(-1:1:1:-1)\\
(R_4A_2,R_4A_2):&-1&-\frac{1}{2}&(1:0:0:0)&(0:0:1:1)\\
x_1=0, x_2=x_3&-\frac{1}{3}&&(\pm 1:0:1:1)&\\
(S,S):&-1&-\frac{1}{3}&(1:0:0:0)&(0:1:1:1)\\
x_1=x_2=x_3&-\frac{1}{4}&&(\pm 1:1:1:1)&\\
\end{array}
\end{eqnarray*}

where we do not write the matrix $(R_4,R_4)$, since the fix lines of the matrices in its conjugacy class are the same as those of the matrices in the conjugacy class of $(A_2,A_2)$. From this table follows that:\\
\begin{itemize}
\item the fix lines above meet different surfaces, hence the conjugacy classes of these matrices are in fact, all distinct (cf. also \cite{sa} (7.3)).\\
\item The fix lines contain the maximal number possible (four), of nodes (this shows {\it 2.})
\end{itemize}



About the number of fix lines: observe that the conjugacy classes of $(S,S)$ has order $32$ and contain the elements $(S^2,S^2)$ hence we have $\frac{32}{2}=16$ distinct fix lines. Finally the conjugacy class of $(A_2,A_2)$ contains six elements, $(A_2,A_1)$ contains twelve elements and $(R_4A_2,R_4A_2)$ contains $24$ elements. Since each element in these conjugacy classes has two fix lines and the elements with minus sign are in the same conjugacy class, the number of fix lines is the same as the number of matrices.
\eprf
 We know that the singular points form $[3,4,3]$-orbits, but we can now show something more. \\
\begin{lemma} the nodes on each singular surface form one $(OO)'$-orbit:

\begin{eqnarray*}
\begin{array}{l|l|l|l|l}
\lambda&-1&-\frac{1}{2}&-\frac{1}{3}&-\frac{1}{4}\\
\hline
{\it orbit}&4&12&16&8\\
\hline
{\it fix~group~mod.}\pm 1&S_4&D_4&D_3&A_4\\      
\hline
{\it order} &24&8&6&12\\
\end{array}
\end{eqnarray*}
\end{lemma}

{\it Proof.} The situation is easy for $\lambda=-1$. In fact the matrix $C$$\in$$[3,4,3]$ let each singular point fix, hence the group $(OO)'$ must permute them. Consider $\lambda=-\frac{1}{2}$, $-\frac{1}{3}$ resp. -$\frac{1}{4}$ and assume that the orbit's length of  singular points is less or equal then $12$, $16$, resp. $8$. Then the fix group in $(OO)'$ mod. $\pm 1$ has order bigger or equal then $8$, $6$ resp. $12$. Checking in the table given in the previous page we see that in fact the only possibility is to have equality.
\eprf

Put now $N_0=$number of nodes on a surface in the pencil, $N_1=$number of fix lines of matrices in the same conjugacy class, $n_0$=number of nodes on a line, $n_1$=number of line through a point. Knowing the fix groups of the singular points and using the formula
\begin{eqnarray}\label{conf}
N_0\cdot n_1=N_1\cdot n_0
\end{eqnarray}
for a configuration of lines and points (cf.\ \cite{sa} section 11.), we can write the table:

\begin{eqnarray*}
\begin{array}{c|c|c|c|c|c|c}
{\rm Repr.~of~the~conj.~class}&\multicolumn{3}{c|}{{\rm value~of }~\lambda}&\multicolumn{3}{c}{{\rm Configuration}}\\
\hline
(A_2,A_2)&-1&-\frac{1}{2}&&(4_3,6_2 )&(12_1,6_2)&\\
(A_2,A_1)&-\frac{1}{2}&-\frac{1}{4}&&(12_2,12_2)&(8_3,12_2)&\\
(R_4A_2,R_4A_2)&-1&-\frac{1}{2}&-\frac{1}{3}&(4_6,24_1)&(12_2,24_1)&(16_3,24_2)\\
(S,S)&-1&-\frac{1}{3}&-\frac{1}{4}&(4_4,16_1)&(16_1,16_1)&(8_4,16_2)\\
\end{array}
\end{eqnarray*}
where in the second column we mean the nodes of the surfaces with the given $\lambda$.



\subsection{The pencil of $IO$-invariant 12-ics}
We take the generators $q^6$ and $S_{IO}$. First we compute the base locus. We use the terminology of \cite{sa} definitions  (5.1),  (5.2) and the facts (5.3) and (5.4). Moreover we consider the groups $(1,O)$, resp. $(I,1)$, which operate on the two rulings of the quadric $q$. To convenience of the reader we recall the table on the length of the orbits under the action of the octahedral group $O$, and of the icosahedral group $I$:   

\begin{center}

\renewcommand{\arraystretch}{1.3}

\begin{eqnarray*}
\begin{array}{c|c}
{\rm octahedron} &{\rm icosahedron}\\
\hline
24,~12,~8,~6&60,~30,~20,~12\\
\end{array}
\end{eqnarray*}

\end{center}

\renewcommand{\arraystretch}{1.0}

\noindent
Now we show:\\
\begin{lemma} The variety $q\cap S_{IO}$ consists of an $(I,1)$-orbit of twelve lines and of an $(1,O)$-orbit of six singular lines of $S_{IO}$.\end{lemma}
{\it Proof}. We have deg$(q\cap S_{IO})$=24. By the table above it can only have bi-degree $(12,12)$ or $(12,6)$. An argumentation as in \cite{sa} (5.5),  (b), shows that $q\cap S_{IO}$  splits into the union of lines of the two rulings of $q$. More precisely it contains the $(I,1)$-orbit of twelve lines and the $(1,O)$-orbit of twelve or of six lines. In the last case the lines have multiplicity two in the intersection. 
So take the fix line $(\lambda:\mu:i\lambda:i\mu)$, $(\lambda:\mu)$$\in$$\piu$, of the matrix $(1,A_2)$ in the orbit of length six. A direct computation (with MAPLE) shows that it is singular on $S_{IO}$, so we are done.
\eprf

\begin{lemma}\label{Li} Let $p$ be a singular point on a surface of the pencil $S_{IO}+\lambda q^6$ and assume that $p$ is on the quadric $q$. Then $p\in L_i$, $i=1,\cdots, 6$, where the $L_i$'s denote the singular lines of the base locus.\end{lemma}
{\it Proof.} As in the proof of  (6.1) in \cite{sa}, if $p$ is singular on a surface in the pencil and $p\in q$, then $p$ is a singular point of $q\cap S_{IO}$. Hence $p\in L_i$.
\eprf
By this fact and Bertini's theorem follows that:\\
\begin{lemma}  1. the general surface in the pencil is smooth away from the lines $L_1,\ldots,L_6$.\\
Moreover we have:\\
2. the singular surfaces have only isolated singularities away from the singular lines $L_1,\ldots,L_6$. In particular they are irreducible and reduced.\end{lemma}
{\it Proof of 2.} 
The surface $q^6$ is the only not reduced surface in the pencil. Indeed, assume that there is another not reduced surface $\{f_{\alpha}^m=0\}$, $\alpha\cdot m=12$, and deg$f_{\alpha}=\alpha$. Then the twelve lines of the base locus in the intersection of $f_{\alpha}^m$ and $q$ would be singular too, which is not the case. Assume now that there is a surface $S$ in the pencil which contain a singular curve $C$. The latter meets the quadric in at least one point, which by lemma \ref{Li} is on some $L_i$. This point is a ``pinch point'' of $S$ (for the definition of ``pinch point'' cf.\ \cite{g-h} p.\ 617). We compute explicitly the pinch points along a line $L_i$. Because of symmetry reason it is enough to do it just for one line. We take $L_1$:$\{x_2=ix_0,x_3=ix_1\}$. First we consider the transformation $x_0\mapsto y_0$, $x_1\mapsto y_1$, $x_2\mapsto y_2+i\cdot y_0$, $x_3\mapsto y_3+i\cdot y_1$, which maps $L_1$ to $\{y_2=y_3=0\}$. We call the transform $L_1$ again.\\
Let $F:=F(y_0,y_1,y_2,y_3,\lambda)$ be the equation of the pencil after the transformation. We can write the Taylor expansion of $F$ along $L_1$
\begin{eqnarray*}
F=A(y_0,y_1,\lambda)\cdot y_2^2+B(y_0,y_1,\lambda)\cdot y_2\cdot y_3+C(y_0,y_1,\lambda)\cdot y_3^2\\
+{\rm terms~of~order}\geq 3,
\end{eqnarray*}
the pinch points are solution of 
\begin{eqnarray*}
{\rm det}\left( \begin{array} {cc}
A& \frac{1}{2}B\\
\frac{1}{2}B& C
\end{array} \right)=0
\end{eqnarray*}
which splits into the product
\begin{eqnarray*}
\begin{array}{l}
(x_1^4+2x_1^3x_0+2x_0^2x_1^2-2x_0^3x_1+x_0^4)(x_1^4-2x_1^3x_0+2x_0^2x_1^2+2x_0^3x_1+x_0^4)\\(x_0^2+3x_0x_1+\sqrt{5}x_0x_1-x_1^2)(x_0^2-3x_0x_1-\sqrt{5}x_0x_1-x_1^2)\\
(3x_0^4-2\sqrt{5}x_0^2x_1^2+3x_1^4)(x_0^4+16x_0^2x_1^2-6\sqrt{5}x_0^2x_1^2+x_1^4)=0.\\
\end{array}
\end{eqnarray*}
this has twenty distinct solutions (independent from $\lambda$), therefore we have ``simple'' pinch points on $L_i$, hence no singular curves.

\eprf

Observe that (6.3) of \cite{sa} holds for this pencil too. In particular, we get the estimate for the number of isolated singularities as follows (we use the same notation as there): deg$(S\cap \partial_i S)=12\cdot 11$ and $S\cdot \partial_i S=C+2L_1+\ldots+2L_6$ hence deg$C=12\cdot 11-12$. Since $C$ cannot be singular there is a $j=0,1,2,3$ s.\ t.\ $C$$\not\subseteq$$\partial_j S$. Hence deg$(C\cap \partial_j S)=(12\cdot 11-12)(12-1)=120\cdot 11$. Since the singular points are computed two times in the intersection their number is $\leq \frac{120\cdot 11}{2}$$=660$.\\
Now  (7.1) of \cite{sa} holds with mult$_{z_{ij}}(L\cdot S)=2$. In fact a line like $L$ there can be only a fix line of a matrix in the conjugacy class of $(A_2,A_2)$ and it meets the base locus at some line $L_i$. Similar assertion as (7.2),  (7.3),  (7.5) and the first part of section 8. hold here too. About the last part of section  8.: in the case of the map $\bar f$, the points $z_{13}$ and $z_{24}$ have order $\frac{12}{2}-2$, hence the degree of the ramification locus is $2\cdot 12 -10$ and without the two points on the quadric, it is $12-4=8$. In particular we have again an estimation for the number of singular points on the fix lines, i. e. this is $\leq$$12$ if they do not meet the base locus, $\leq 8$ otherwise.\\
{\it The singular surfaces.} 
We proceed by a direct computation using the lines of fix points. 
The matrices of $IO$ with fix lines containing singular points are in the conjugacy class of $(A_2,A_2)$, $(S,S)$ or $(A_2,R_4A_2)$. The total number of distinct fix lines is 90, 80,  resp. 180. In the table below we give the singular surfaces and the singular points on the fix lines. We choose a representative in each conjugacy class and a fix line of it, moreover we put $a':=\sqrt{2}-1$, $a:=\sqrt{2}+1$,$\gamma_1=2\sqrt{2}-\sqrt{5}$, $\gamma_1'=2\sqrt{2}+\sqrt{5}$, $\gamma_2:=2+\sqrt{2}+\sqrt{5}+\sqrt{10}$, $\gamma_2':=-2+\sqrt{2}-\sqrt{5}+\sqrt{10}$,  $\gamma_3:=-2+\sqrt{2}+\sqrt{5}-\sqrt{10}$, $\gamma_3':=2+\sqrt{2}-\sqrt{5}-\sqrt{10}$, $\alpha_1:=\frac{1}{2}(-1+\sqrt{2}-\sqrt{10})$, $\alpha_1':=\frac{1}{2}(1+\sqrt{2}-\sqrt{10})$, $\beta_1:=\frac{1}{2}(4-3\sqrt{2}-2\sqrt{5}+\sqrt{10})$, $\alpha_2:=\frac{1}{2}(1+\sqrt{2}+\sqrt{10})$, $\alpha_2':=\frac{1}{2}(-1+\sqrt{2}+\sqrt{10})$, $\alpha_3:=-2+\frac{3}{2}\sqrt{2}-\sqrt{5}+\frac{1}{2}\sqrt{10}$, $\alpha_3':=2+\frac{3}{2}\sqrt{2}+\sqrt{5}+\frac{1}{2}\sqrt{10}$, $\beta_2:=\frac{1}{2}(-10+7\sqrt{2}-4\sqrt{5}+3\sqrt{10})$, $\beta_3:=\frac{1}{2}(\sqrt{2}-2\sqrt{5}+\sqrt{10})$, $\beta_4:=\frac{1}{2}(2+\sqrt{2}+\sqrt{10})$, $c_1:=-\frac{74}{972}+\frac{4}{243}\sqrt{10}$, $c_2:=-\frac{74}{972}-\frac{4}{243}\sqrt{10}$:


\vspace*{0.3cm}

\begin{eqnarray*}
\begin{array}{c|c|c|c|c}
{\rm matrix~and~fix~line}&\multicolumn{2}{c|}{{\rm value~of} ~\lambda}&\multicolumn{2}{c}{{\rm int.~points}}\\
\hline
(A_2,A_2): &0&-\frac{1}{8}&(1:0:0:0)&(\pm a':0:1:0)\\
x_1=x_3=0&&&(0:0:1:0)&(1:0:\pm a':0)\\
&&&(\pm 1:0:1:0)&\\
(A_2,R_4A_2):&-\frac{1}{8}&&(-a:1:0:0)&\\
x_0=-ax_1, x_3=a'x_2&&&(0:0:1:a')&\\
&&&(1:\pm a':1:a')&\\
&c_1&c_2&(\alpha_1:-\beta_1:1:a')&(\alpha_1':-\beta_3:1:a')\\
&&&(-\alpha_2:\beta_1:1:a')&(-\alpha_2':\beta_3:1:a')\\
&&&(\alpha_3:\beta_2:1:a')&(\alpha_3':\beta_4:1:a')\\
&&&(-\alpha_3:-\beta_2:1:a')&(-\alpha_3':-\beta_4:1:a')\\
(S,S):&0&&(1:0:0:0)&\\
x_1=x_2=x_3&&&(0:1:1:1)&\\
&&&(-1:1:1:1)&\\
&&&(-\sqrt{5}:1:1:1)&\\
&&&(\sqrt{5}\pm 2:1:1:1)&\\
&c_1&c_2&(\gamma_1:1:1:1)&(\gamma_1':1:1:1)\\
&&&(\gamma_2:1:1:1)&(\gamma_2':1:1:1)\\
&&&(\gamma_3:1:1:1)&(\gamma_3':1:1:1)\\
\end{array}
\end{eqnarray*} 
Here the fix line of $(A_2,A_2)$ contain two points of the base locus. Observe that the number of singular points is maximal on the lines, hence the four given $\lambda$ are the only values corresponding to singular surfaces.\\
\begin{prop} In the pencil $S_{IO}+\lambda q^6$ we have the following $IO$-orbits of nodes

\begin{center}
\begin{tabular}{l|*{4}{l}}
$\lambda$&$0$&$c_1$&$-\frac{1}{8}$& $c_2$\\
\hline
orbit&$120$&$240$&$360$&$240$\\
\hline
fix group mod. $\pm 1$& $A_4$&$D_3$&$\mathbb{Z}_2\times\mathbb{Z}_2$&$D_3$\\   
\hline
order&$12$&$6$&$4$&$6$\\
\end{tabular}
\end{center}
 \end{prop}
\bprf 
 We analyze case for case the four different surface. A direct computation  as in section  4.2, then shows that the singular points are nodes.\\
I. \underline{$\lambda=0$}. Observe that $(1:0:0:0)\in$$S_{IO}$ is contained in three fix lines of elements in the conjugacy class of $(A_2,A_2)$. Moreover \\

(*) there is a matrix of order four $(1,\gamma)$$\in IO$, with $\gamma^2=A_2$, \\

Hence $(1,\gamma)$ and $(A_2,A_2)$ commute and the four points on the fix line of $(A_2,A_2)$ form one orbit under the action of $(1,\gamma)$. By this fact follows that we have a configuration of lines and points. By the formula (\ref{conf}), we get $N_0\cdot 3=90\cdot 4$, hence $N_0=120$. The point $(1:0:0:0)$ is on the fix line of $(S,S)$ too, more precisely it is contained on four such lines. The formula $120\cdot 4=80\cdot n_0$ shows $n_0=6$ so the six points on the fix line are in the orbit of length 120 too. In conclusion we have just one $IO$-orbit of singular double points.\\
II. \underline{$\lambda=-\frac{1}{16}$}. The four points on the fix line of $(A_2,A_2)$ form one orbit by (*). The formula (\ref{conf}) gives $N_0\cdot n_1=360$. If $n_1\geq 2$ then $N_0\leq 180$ and the fix group in $IO$ mod. $\pm 1$ has order $\geq 8$. Checking in the table on page 24 this is not possible. Hence $N_0=360$ and $n_1=1$. The fix group has order four and it is isomorphic with  $\mathbb{Z}_2\times\mathbb{Z}_2$. Each singular point is contained in two fix lines of matrices in the conjugacy class of $(A_2,R_4A_2)$. Now the formula (\ref{conf})) with $N_0=360$, $n_1=2$ and $N_1=180$ gives $n_0=4$. So we conclude that the singular points on these fix lines are in the orbit of length $360$ too.\\  
III. \underline{$\lambda=c_i,~i=1,2$}. The three points on the fix line of $(S,S)$ form one orbit under the action of $(S,1)$. The formula (\ref{conf}) gives $N_0\cdot n_1=240$. If $n_1\geq 2$ then $N_0\leq 120$ and the fix group mod. $\pm 1$ has order $\geq 12$, which is not possible (check in the table on page 24 again). Hence $N_0=240$ and $n_1=1$. The fix group has order six and it is isomorphic with $D_3$, hence three fix lines of matrices in $(A_2, R_4A_2)$ contain a singular points. The formula (\ref{conf}) in this case gives $n_0=4$, which shows that the singular points on these fix lines are in the previous orbit too.

\eprf

We give now the table of the configurations of lines and points:

\vspace*{0.3cm}

\renewcommand{\arraystretch}{1.3}

\begin{eqnarray*}
\begin{array}{c|c|c|c|c}
{\rm Repr.~of~the~ conj.~class}&\multicolumn{2}{c|}{{\rm value~of }~\lambda}&\multicolumn{2}{c}{{\rm Configuration}}\\
\hline
(A_2,A_2)&0&-\frac{1}{8}&(120_3,90_4)&(360_1,90_4)\\
(A_2,R_4A_2)&-\frac{1}{8}&c_i&(360_2,180_4)&(240_3,180_4)\\
(S,S)&0&c_i&(120_4,80_6)&(240_1,80_3)
\end{array}
\end{eqnarray*}

\subsection{The pencil of $D_nD_m$-invariant quartics}
We take the generators $q^2$ and $f_2+f_3=(x_0-ix_1)(x_0+ix_1)(x_2-ix_3)(x_2+ix_3)$. The latter is the union of four complex planes meeting each other at the four complex lines $\{x_0\pm ix_1=0,~x_2\pm ix_3=0\}$ on $q$. 
The pencil contain the multiple quadric $q'^2$$:=q^2-2(f_2+f_3)$$=(x_0^2+x_1^2-x_2^2-x_3^2)^2$ too. Hence with this new generator the equation of the pencil becomes
\begin{eqnarray*}
q'^2+\lambda q^2=(q'+i\sqrt{\lambda}q)(q'-i\sqrt{\lambda}q).
\end{eqnarray*}
A surface in the pencil is the union of two quadrics which meet each other at the four lines given above. An easy computation shows that none of the surface in the pencil contains isolated singularities.

\input{artjul01_22}

\input{i_o}



\addcontentsline{toc}{section}{  \hspace{0.5ex} References}

%% file: artjul01_22.tex
\subsection{The pencil of $OD_n$-invariant octics ($n\geq 6$) and of $ID_n$-invariant 12-ics ($n\geq 8$).}
We take generators $q^4$, $q^6$ and $P_8=\prod_{i=0}^{3}K_i$, $P_{12}=\prod_{i=0}^{5}K_i'$ (cf. {\bf (3.1.2)}) and we denote the pencils by $\Pi_m(\lambda)$$:=P_m+\lambda q^{\frac{m}{2}}$, $m=8,12$. We have\\
\begin{lemma} The base locus of the pencil is $\frac{m}{2}$-times the intersection $q\cap P_m$, $m=8,12$. Where (as sets)
\begin{center}
\begin{eqnarray*}
\begin{array}{lll}
q\cap P_8&=& {\rm fix~lines~of~}(1,R_n),\\ && {\rm orbit~of~eight~lines~under}~(O,1),\\
q\cap P_{12}&=&{\rm fix~lines~of~}(1,R_n), \\ && {\rm orbit~of~twelve~lines~under}~(I,1).\\
\end{array}
\end{eqnarray*}
\end{center}
\end{lemma}
\bprf
The element $(1,R_n)$ has fix lines $L_1:\{x_0=-ix_1,~x_2=ix_3\}$, $L_2:\{x_0=ix_1,~x_2=-ix_3\}$ for each $n$. These are contained in each quadric $K_s$, $s=0,\ldots,3$, resp. $K_t'$, $t=0,\ldots,5$ hence they have at least multiplicity four, resp. six in the intersection $q\cap P_m$, $m=8,12$. By this fact and an argumentation as in \cite{sa}, (5.5), it follows that the lines of the length $m=8,12$ orbit under $(O,1)$, resp. $(I,1)$ are in the base locus too. Moreover the previous multiplicity of intersection are exactly four, resp. six. \eprf

\begin{lemma} The pencil does not contain surfaces with isolated singular points. \end{lemma}
\bprf
The groups have order $96 n$, resp.\ $240n$, since an isolated  singular point has orbit of finite length, it is fixed by $(1,R_n)$. This shows that it is  on the lines of the base locus, hence not isolated. A contradiction to the assumption.\eprf
\begin{prop} In the pencil we have the following singular surfaces, with one orbit of double lines which are fix lines for the elements of some conjugacy class:
\begin{eqnarray*}
\begin{array}{c|cccc}
{\rm value~of}~\lambda& OD_n:& -1&\frac{1}{3}&0\\
&ID_n:&0&&\\
\hline
{\rm orbit}&OD_n:&6&8&12\\
&ID_n:&30&&\\
\hline
{\rm fix~lines~in }&&(A_3,R_n^{\frac{n}{2}})&(S,R_6^2)&(R,R_n^{\frac{n}{2}})\\
{\rm the~conj.~class~of}&&&&\\
\end{array}
\end{eqnarray*}
\end{prop}
\bprf
It is a direct computation, using the equations $x_0=x_2,x_1=-x_3$; $x_0=-(1+\sqrt{3})x_3+x_1, x_2=x_3+(1-\sqrt{3})x_1$ and $x_2=(1-\sqrt{2})x_0,x_3=(\sqrt{2}-1)x_1$ of a fix line of ($A_3,R_n^{\frac{n}{2}}$), ($S,R_6^2$), resp.\ ($R,R_n^{\frac{n}{2}} $). \eprf
\section{Final Remarks}
1) By \cite{cox} p.\ 292 there are sixteen regular polytopes $\{p,q,r\}$ in four dimensions. These polytopes correspond to four distinct symmetry groups $[p,q,r]$ listed in the table below. Three of these symmetry groups can be obtained from groups which we described above by adding extra generators $C$ and $C'$, where $C$ denotes the matrix given in section \ref{invariants}, $C'$ is the matrix of \cite{sa} p.\ 433.
\begin{eqnarray*}
\begin{array}{c|c|c|c|c}
{\rm symmetry~groups}&[3,3,3]&[3,3,4]&[3,4,3]&[3,3,5]\\
\hline
{\rm our~groups}&&(OO)'&TT&II\\
\hline
{\rm extra~generators}&&C&C,C'&C\\
\end{array}
\end{eqnarray*}
The group $[3,3,3]\cong Sym(5)$ does not contain the Heisenberg group $H$, in fact it has some invariant polynomials of odd degree (cf.\ \cite{cox0} p.\ 780). In these notes we complete the description of the $G$-invariant pencils of surfaces whenever $G$ is the symmetry group of a regular  four dimensional polytope and $G$ contains $H$.\\
2) In remark 3.3, with $n=3$, we give the degree four invariant polynomials of the {\it bi-polyhedral dihedral group} $D_3D_3$$\subseteq$$SO(4)$. By Mukai \cite{mukai} the quotient $\pitr$$/D_3D_3$ is isomorphic with the Satake compactification of the moduli space of abelian surfaces with $(1,2)$-polarization, hence this invariant polynomials should be related to modular forms.\\
3) In \cite{basa} the quotients $X/G$, $G=TT,OO,II$ are described. It would be  interesting to examine the quotients in the remaining cases.

%% file: i_o.tex
\section{Computer Picture}
We exhibit a computer picture of the $I\times O$ symmetric surface of degree 12 with 360 nodes. This has been realized with the program SURF written by S.\ Endra\ss.
%
\hoffset=1.5cm
\voffset=1.0cm
\begin{center}
\hspace*{-4cm}
\epsfig{file=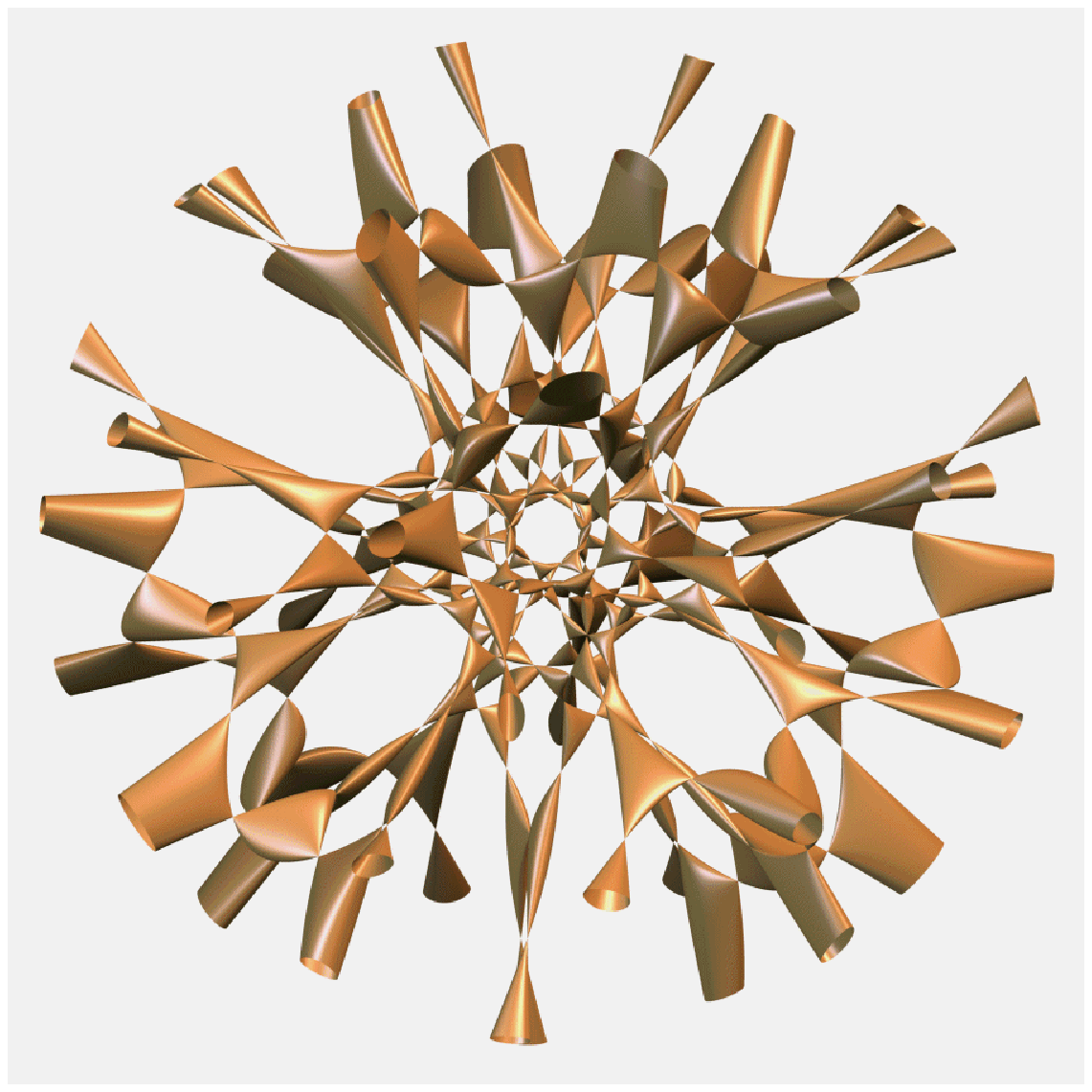}
\hspace*{-4cm}\ \\[3ex]
 {\large $I\times O$-symmetric surface of degree 12}\\[0.5cm]
{\large with 360 nodes}

%
\end{center}
%